%% file: design_space.tex
\documentclass{article}

\usepackage{arxiv}

\usepackage[utf8]{inputenc} 
\usepackage[T1]{fontenc}    
\usepackage[hidelinks]{hyperref}       
\usepackage{url}            
\usepackage{booktabs}       
\usepackage{amsfonts}       
\usepackage{nicefrac}       
\usepackage{microtype}      
\usepackage{lipsum}		
\usepackage{graphicx}
\usepackage[numbers]{natbib}
\usepackage{doi}
\usepackage{caption}
\usepackage{subcaption}
\usepackage{float}
\usepackage{amsmath}

\title{A Method for Finding a Design Space as Linear Combinations of Parameter Ranges for Biopharmaceutical Control Strategies}

\author{ \href{https://orcid.org/0000-0002-1684-5150}{\includegraphics[scale=0.06]{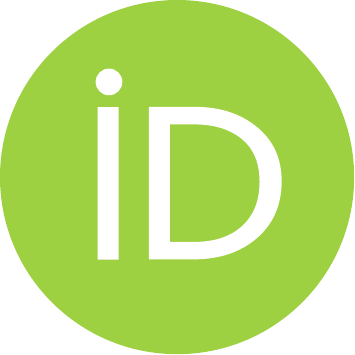}\hspace{1mm}Thomas Oberleitner}\\
	Competence Center CHASE GmbH\\
    Ghegastraße 3, Top 3.2, 1030 Vienna, Austria
	\And
	Thomas Zahel\\
	Körber Pharma Austria GmbH, PAS-X Savvy\\
	Mariahilferstr. 88A/1/9, 1070 Vienna, Austria\\
    \And
	\href{https://orcid.org/0000-0003-2314-1458}{\includegraphics[scale=0.06]{orcid.pdf}\hspace{1mm}Christoph Herwig}\\
	Research Area Biochemcial Engineering\\
	Vienna University of Technology\\
	Getreidemarkt 9, 1060 Vienna, Austria\\
	\texttt{christoph.herwig@tuwien.ac.at}\\
}

\renewcommand{\shorttitle}{Design Spaces as Linear Combinations of Parameter Ranges}

\begin{document}
\maketitle

\begin{abstract}
According to ICH Q8 guidelines, the biopharmaceutical manufacturer submits a design space (DS) definition as part of the regulatory approval application, in which case process parameter (PP) deviations within this space are not considered a change and do not trigger a regulatory post approval procedure. A DS can be described by non-linear PP ranges, i.e., the range of one PP conditioned on specific values of another. However, independent PP ranges (linear combinations) are often preferred in biopharmaceutical manufacturing due to their operation simplicity. While some statistical software supports the calculation of a DS comprised of linear combinations, such methods are generally based on discretizing the parameter space - an approach that scales poorly as the number of PPs increases. Here, we introduce a novel method for finding linear PP combinations using a numeric optimizer to calculate the largest design space within the parameter space that results in critical quality attribute (CQA) boundaries within acceptance criteria, predicted by a regression model. A precomputed approximation of tolerance intervals is used in inequality constraints to facilitate fast evaluations of this boundary using a single matrix multiplication. Correctness of the method was validated against different ground truths with known design spaces. Compared to state-of-the-art, grid-based approaches, the optimizer-based procedure is more accurate, generally yields a larger DS and enables the calculation in higher dimensions. Furthermore, a proposed weighting scheme can be used to favor certain PPs over others and therefore enabling a more dynamic approach to DS definition and exploration. The increased PP ranges of the larger DS provide greater operational flexibility for biopharmaceutical manufacturers.
\end{abstract}


\input{sections/01_introduction}
\input{sections/02_methods}
\input{sections/03_contribution}
\input{sections/04_sim_results}
\input{sections/05_discussion}
\input{sections/06_conclusion}
\input{sections/07_acknowledgements}

\bibliographystyle{unsrtnat}
\bibliography{design_space}  

\end{document}

%% file: sections/01_introduction.tex
\section{Introduction}
\begin{figure}[ht]
    \begin{subfigure}[t]{0.50\textwidth}
        \centering
        \includegraphics[]{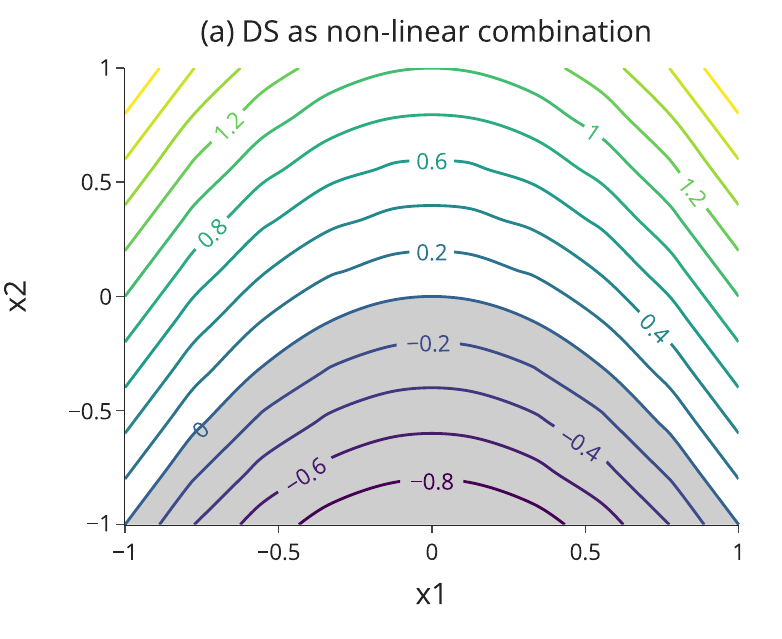}
    \end{subfigure}
    \hfill
    \begin{subfigure}[t]{0.50\textwidth}
        \centering
        \includegraphics[]{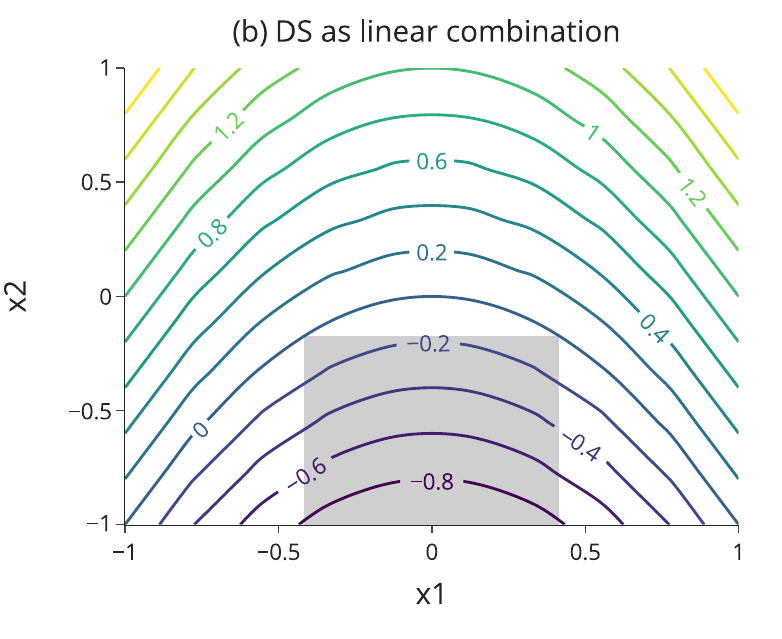}
    \end{subfigure}

\caption{The design space for the function $f(x_1,x_2 )=x_1^2+x_2$, where $f(x_1,x_2 )\leq0$, shown in the contour plot as a non-linear (a) and linear (b) combination of input parameters $x_1$ and $x_2$.}
\label{fig:ds-lincomb}
\end{figure}

The ICH Q8 guideline for pharmaceutical development defines the design space (DS) as “the multidimensional combination and interaction of input variables (e.g., material attributes) and process parameters that have been demonstrated to provide assurance of quality” \cite{ich2017q8r2}. The process parameters (PP) described here are generally identified in the risk assessment or process development phases and are considered critical when sufficient evidence was found that they affect the output of a unit operation, i.e., a critical quality attribute (CQA). A design space is comprised of the ranges of these process parameters that result in CQA values within acceptable limits. For the biopharmaceutical manufacturer a DS definition can be submitted as part of the regulatory approval application, in which case PP deviations within this space are not considered a change and therefore do not trigger a regulatory post approval procedure. For operators, the DS constitutes a valuable guideline document for controlling a process. While ICH Q8 does not recommend a specific form or method for calculating a DS, it provides examples for how to present the DS as non-linear and linear combinations of parameter ranges in the form of contour plots (appendix 2c in guide). Non-linear combinations describe the DS as a set of rules, or parameter ranges conditioned on other parameters, e.g., “PP1 is allowed to move between -1 and 1 if PP2 is lower than 0.5”. Linear combinations of parameter ranges on the other hand are independent of each other. While the former description generally represents a larger space to operate in and methods for computing it can be found several publications \cite{kim2002design,kusumo2020nested}, the latter might be preferred due to its operational simplicity and is the subject of this contribution. Figure \ref{fig:ds-lincomb} shows the different types of design space graphically.

In these examples, the definition of a design space does not incorporate any measure of statistical uncertainty, i.e., the contour shown in Figure \ref{fig:ds-lincomb} directly represents the predicted mean CQA values from the model. In the context of biopharmaceutical process validation, we suggest a more conservative approach. To accurately quantify uncertainty inherent in the regression model due to analytical and process variability, we replace predicted CQA values with the upper and lower boundary of a tolerance interval (TI) that incorporates nominal levels of both confidence and coverage. The statistical relevance of tolerance intervals in biopharmaceutical control strategies has been discussed earlier \cite{oberleitner2023incorporating}. These boundaries are then used to validate acceptance limits in a conservative manner.

As the visual definition of a DS is only feasible in a low-dimensional parameter space, e.g., the bivariate contour plots in Figure \ref{fig:ds-lincomb}, computational methods are required to find the exact parameter ranges for typical biopharmaceutical models containing 10 parameters or more. To the knowledge of the authors, this is currently not possible in state-of-the-art statistical software used in industry because the common approach is to discretize the parameter space and to evaluate points on the resulting grid. This method suffers from the curse of dimensionality as the parameter space grows exponentially over the number of parameters in the model. Segmenting the range of 10 parameters into 10 parts would result in a parameter space of $10^10$ points, for which all possible combinations would need to be evaluated to find the largest possible DS – an exceptionally computationally expensive task. As a result, this makes it practically impossible to calculate design spaces for more than 6-7 PPs, however, in typical scenarios of biopharmaceutical process development and characterization this number is exceeded. Furthermore, the design space found by grid-based methods is generally not the one with the largest possible multidimensional volume, or “hypervolume”, because the solution space is limited to discrete points that are distributed over the parameter range. We illustrate these drawbacks in section 4. As a result, currently, the elaborated DSs from large Quality by Design (QbD) projects are not used in the control strategy and the manufacturer goes back to univariate individual controls, which leads to a loss of flexibility and missing the opportunity of process optimization to maneuver in the DS.

In this contribution, we present a novel method for finding a design space comprised of linearly independent parameter ranges while treating CQA predictions conservatively by evaluating the tolerance interval boundaries around them, i.e., checking whether those boundaries fall within acceptance criteria. We propose a method based on numeric optimization rather than the screening of a grid, which results in improved computing times and, in many cases, a larger volume of the DS. An approximation for tolerance intervals is presented as well as a method for incorporating categorical model factors in a computationally efficient manner. Finally, a weighting scheme to favor certain PPs over others facilitates increased flexibility in the calculation of the results. The overall result is a method for DS computation and exploration that is computationally more efficient and provides greater flexibility than the state-of-the-art, both in terms of controlling the DS computation as well as greater flexibility for operators due to increased PP ranges.

%% file: sections/02_methods.tex
\section{Methods}
\subsection{Regression Models in Biopharmaceutical Development and Manufacturing}
Due to their simplicity and optimal statistical properties, regression models are a popular choice in biopharmaceutical development and manufacturing, where they are used to express the relationship between PPs and CQAs. One of its more basic representatives is the ordinary least-squares (OLS) model, which assumes a (curvy-) linear relationship between parameters and response and normally distributed residuals \cite{montgomery2021regression}. The model equation takes the form:

\begin{equation}
y=X\beta+\epsilon
\end{equation}

Where $y$ is the model output (here, the CQA), $X$ an $n*p$ matrix of $n$ observations comprised of $p$ parameter settings (PPs). The formula also includes the vector of model coefficients $\beta$ that assigns each parameter a numeric value (a parameter’s effect) and the residual error $\epsilon \sim N(0, \sigma^2)$. The least-squares fit is obtained by finding an estimator for the model coefficients. In the case of OLS, can be calculated in the following way:

\begin{equation}
\hat{\beta}=(X^TX)^{-1}X^Ty
\end{equation}

Under the assumption that $\epsilon \sim N(0, \sigma^2)$, this definition yields the best linear unbiased estimator (BLUE) for $\beta$. Note that in practice, the assumption of normally distributed residuals might not be satisfied, in which case other types of regression models might be employed where $\hat{\beta}$ is acquired using numerical optimization rather than a closed-form expression \cite{brown2015applied,stroup2012generalized,knudson2016monte}. The method for DS calculation proposed in this article is generally agnostic about the type of model used for representing the relationship between PP and CQA as long as tolerance intervals can be defined for its predictions. For the sake of simplicity and their frequent application in biopharmaceutical process development, we chose OLS models for the simulation studies in 4.

\subsection{Tolerance Intervals}
Tolerance intervals are used to quantify uncertainty of the predicted mean of a model. Its purpose is to estimate the population distribution of the model’s predicted response, given the uncertainty associated with the modelling and sampling process. Using such conservative estimators is especially important in the biopharmaceutical domain where a majority of the future population of runs needs to be within process or specification limits - hence the usage of confidence or prediction intervals is not recommended \cite{oberleitner2023incorporating}. As defined in \cite{degryze2007using}, tolerance intervals estimate the following probability:

\begin{equation}
Pr[Pr[y|x_{n+1} \in I_T(x_{n+1})] \geq P] = 1 - \alpha
\end{equation}

The inner probability expresses whether the model’s response y for a new observation $x_(n+1)$ is contained within the true population distribution $I_T(x_(n+1))$ in a least $100(1-\psi)\%$ of repeated samplings from the reference distribution, whereas the outer probability represents the confidence level of $100(1-\alpha)\%$. Consequently, tolerance intervals contain nominal parameters for both the level of coverage ($\psi$) as well as confidence ($\alpha$).

Depending on the underlying model, the computation of tolerance intervals can be quite complex and might involve numeric optimization, bootstrapping or other computationally expensive methods (for a comprehensive overview, see \cite{krishnamoorthy2009statistical}). This is not the case for OLS models, where tolerance intervals can be calculated using a ratio of critical values of the $\chi^2$ distribution \cite{guenther1972tolerance}:

\begin{equation} \label{eq:ti-approx}
\hat{y} \pm \sigma \sqrt{\frac{(n-p)\chi_{1;\psi}^2(\frac{1}{n_i^*})}{\chi_{n-p;\alpha}^2}}
\end{equation}

The term $\chi_{1;\psi}^2$ in the numerator is the critical value of the $\chi^2$ distribution at probability $\psi$, one degree of freedom and the noncentrality parameter set to $1/n_i^*$ where $n_i^*$ is the vector of “effective number of observations” $n_i^*=\frac{\hat{\sigma}^2}{se(\hat{y}_i)^2}$. The distribution in the denominator is evaluated at probability $\alpha$ and uses $n-p$ degrees of freedom. While not strictly necessary in the case of OLS, we use this definition for illustrating the approximation of intervals in section 3.2.

\subsection{Optimization Algorithms}
The search for the design space with the largest possible volume presented here is formulated as a continuous optimization problem. Such problems find specific values of the parameter vector $x$ that minimizes a function $f(x)$ subject to a set if inequality constraints $g_i (x) \leq 0, i=1,...,m$ where $m$ is the number of constraints. The DS volume problem can be formulated as an optimization objective function and a set of inequality constraints (see section 3 for details). Those relatively minor requirements enable us to choose from a multitude of well researched and widely available optimization algorithms (for benchmarks and an overview see \cite{varelas2019benchmarking}).

COBYLA (Constrained Optimization BY Linear Approximation) was chosen as the main optimization algorithm to minimize the DS volume, as it meets those requirements and, as a gradient-free method, shows reasonable robustness against converging in local minima \cite{powell1994direct}. The algorithm repeatedly evaluates the objective function at the corners of a “simplex”, i.e., at $p+1$ points, $p$ being the number of variables in $x$. One of the most widely used variants of simplex-based optimization algorithms is the one proposed by Nelder and Mead \cite{nelder1965simplex}. The Nelder-Mead algorithm, however, can run into situations where an incorrect minimum of the objective function is found, as shown by MacKinnon \cite{mckinnon1998convergence}. COBYLA expands upon Nelder-Mead by interpolating the vertices of the simplex using linear polynomials and introducing a trust region, which in turn are used to find the next vertex candidate. This new vertex is different from all current vertices in the simplex, thus circumventing the problem of incorrect convergence. The trust region radius $\Delta>0$ represents the boundary for finding a new vertex, i.e., an improved point $\underline{\hat{x}}$ in the vicinity of the current point $\underline{x_0}$ is found by minimizing the objective function $L(\underline{\hat{x}})$ subject to $\|\underline{\hat{x}} - \underline{x_0}\| < \Delta$. The trust region radius is adjusted automatically and associated with a lower boundary $\rho$. This boundary starts with a predefined value $\rho_{start}$ and gets smaller in later iterations to avoid local minima. The optimization parameter $\rho_{start}$ can be set by the user and is relevant for the method proposed in the following sections, as it controls the granularity of the optimization process.

An intermediate step in our proposed approach concerns the search for minima and maxima of the objective function $f(x)$ within the boundaries of a hyperrectangle, i.e., a nested optimization problem with boundaries and no constraints. The L-BFGS-B algorithm was chosen due do its support for simple bounds and performance properties, especially in higher dimensions \cite{varelas2019benchmarking}. This quasi-Newton method uses a limited-memory version of the approximation of the Hessian proposed in the original BFGS algorithm to guide optimization.

We propose an optional, second optimization pass for finetuning results, which can be executed using a COBYLA optimizer with a smaller value of $\rho_{start}$ or SLSQP (sequential least squares programming), a quasi-Newton optimization algorithm \cite{kraft1988software}. This step uses the same constraints as in the main optimization step, except that the true TI calculation is used as opposed to an approximation. The rationale behind this multi-step approach is that the first pass is expected to converge within the vicinity of the global optimum, while the second pass refines results and eliminates potential errors introduced by the TI approximation.

%% file: sections/03_contribution.tex
\section{Finding a Design Space Comprised of Linear Combinations of Parameter Ranges}
\subsection{Overview}
\begin{figure}[ht]
\centering
\includegraphics[]{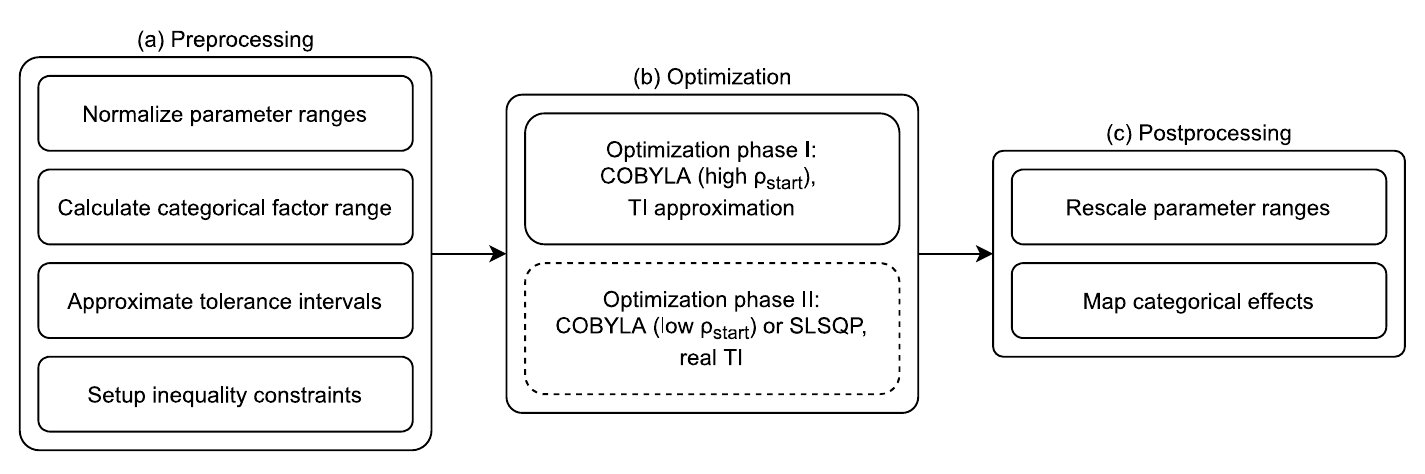}
\caption{Overview of the steps involved in the DS calculation using an optimizer.}
\label{fig:steps}
\end{figure}

The method proposed in this contribution aims to solve the accuracy and computing time problem outlined in the introduction. This is achieved by employing numerical optimizers instead of grid screening methods commonly found in state-of-the-art software. Section 3.3 describes the objective function and inequality constraints used by those optimizers. As the minimization of computing time is of paramount importance, we furthermore introduce a quadratic approximation of tolerance intervals in section 3.2 as well as a fast method for incorporating categorical effects in section 3.4. These steps, as well as the use of an optimization algorithm, require some pre- and post-processing procedures, shown in Figure \ref{fig:steps}.

In the preprocessing step, each dimension of the parameter space is normalized to the range $[-1,1]$ to improve speed and accuracy of the optimization algorithms \cite{horst2013handbook}. Categorical factor ranges are determined, and the tolerance interval approximation is set up. This information and appropriate inequality constraints are then used in the optimization step. The primary optimization phase tries to find the region of the global optimum using COBYLA and large, initial step sizes expressed by the parameter $\rho_{start}$. Optionally, the results from the first phase can be improved in a second phase using either COBLYA with a reduced step size $\rho_{start}$ or, depending on the underlying model and dimensionality, a gradient-based SLSQP optimizer. After a satisfactory DS is found, the results are transformed into their original scale and the continuous boundaries for categorical effects are mapped to valid levels (for details, see next sections).

\subsection{Tolerance Interval Approximation}
\begin{figure}[H]
    \begin{subfigure}[t]{0.50\textwidth}
        \centering
        \includegraphics[width=\textwidth]{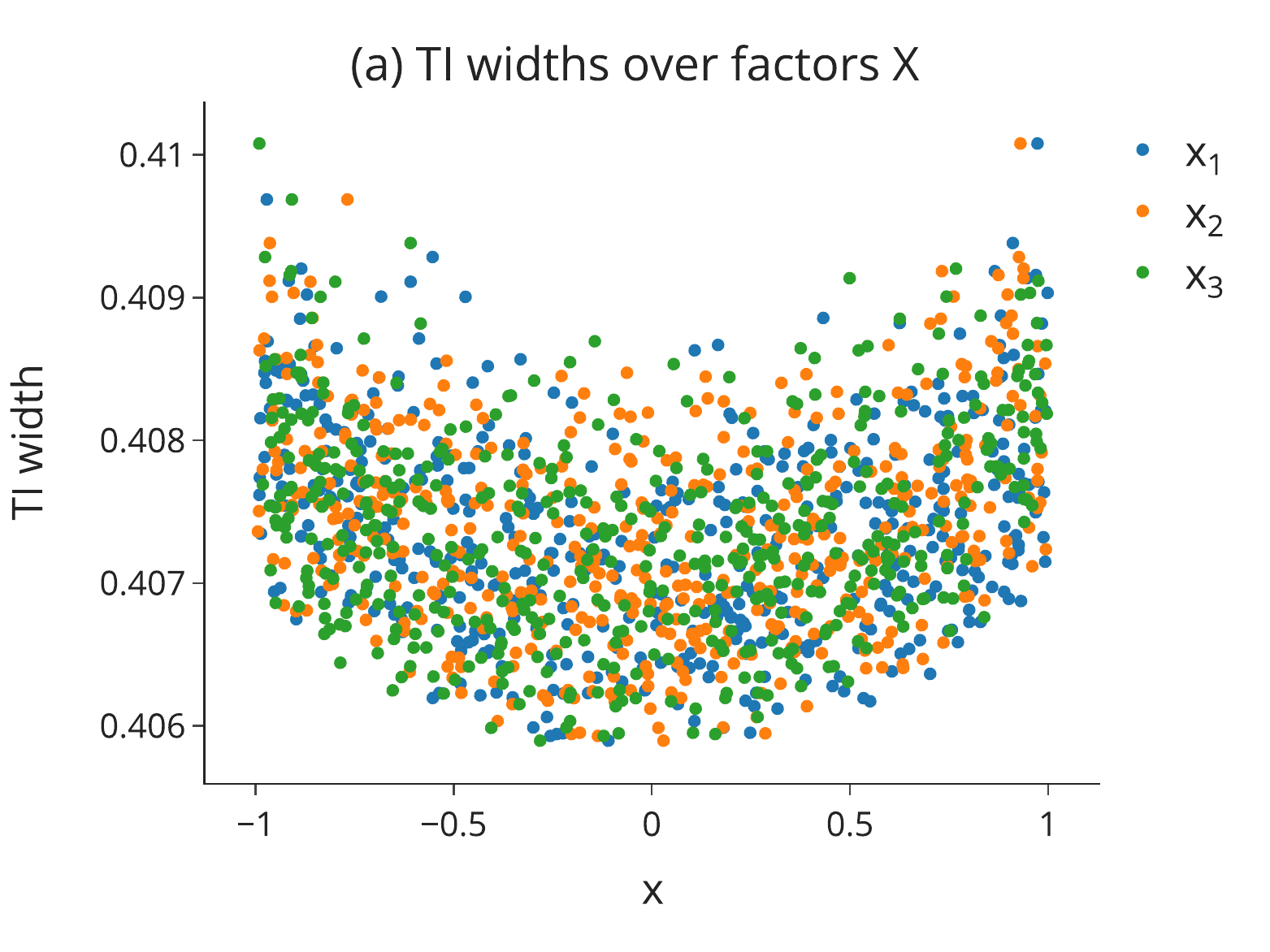}
    \end{subfigure}
    \hfill
    \begin{subfigure}[t]{0.50\textwidth}
        \centering
        \includegraphics[width=\textwidth]{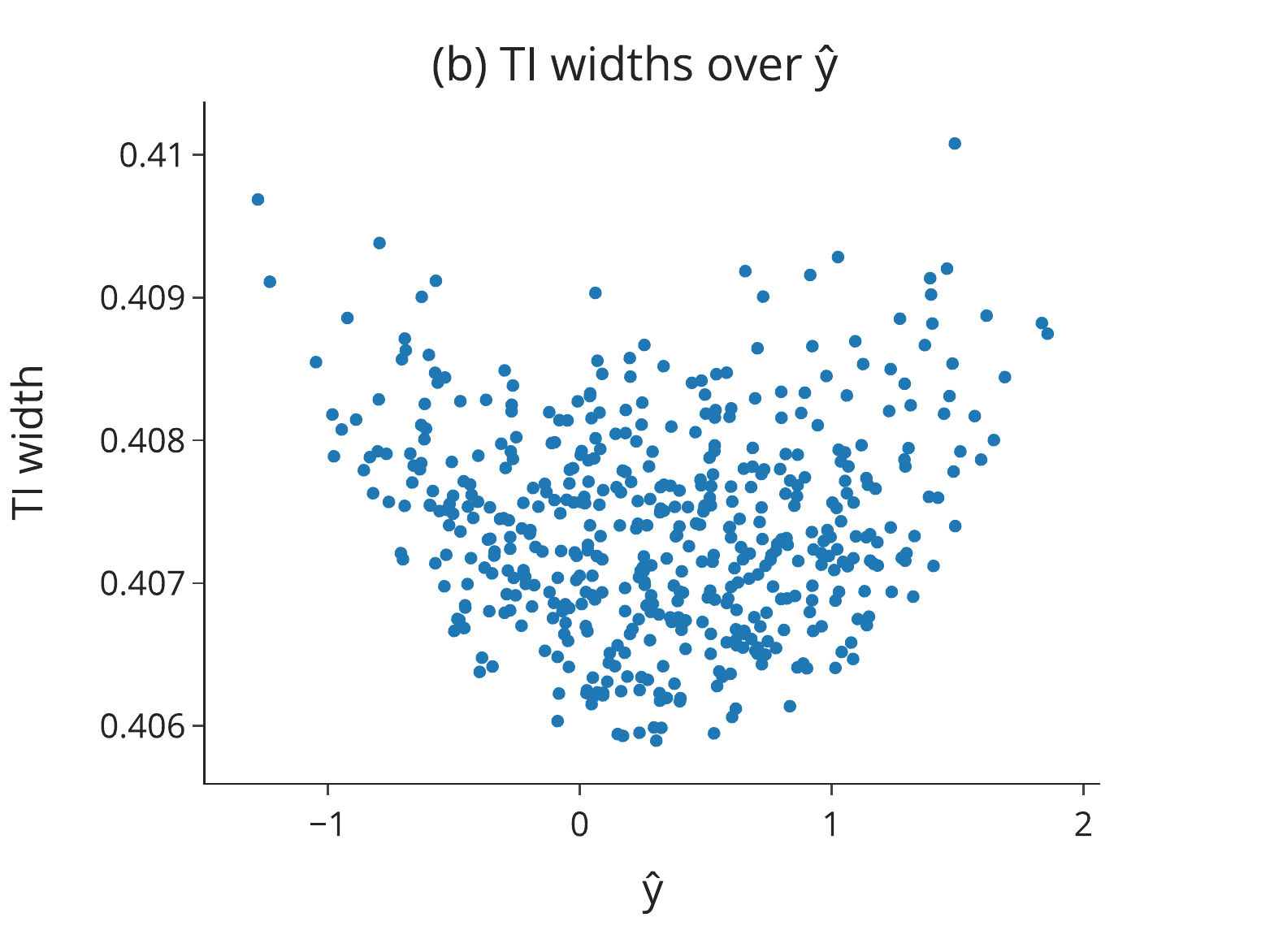}
    \end{subfigure}

\caption{Tolerance interval (TI) widths over ordered (a) factor values and (b) mean predictions.}
\label{fig:ti-shape}
\end{figure}

As mentioned in section 2.2, the calculation of tolerance intervals can be computationally expensive, depending on the type of model used for prediction. Furthermore, the inequality constraints described in section 3.3 evaluate intervals at several points in each iteration of the optimization algorithm, turning them into a potential bottleneck. Therefore, we present an approximation method that turns the calculation of a tolerance interval into a simple vector multiplication that can be carried out in a computationally efficient manner using standard linear algebra libraries. Although not necessarily required in the case of OLS, we use the tolerance interval definition in equation \ref{eq:ti-approx} as an illustrative example. To derive a parsimonious approximation, we exploit a property common to regression models, that is, that the least-squares projection $\hat{y}=X\hat{\beta}$ always goes through the multivariate mean of $X$ and therefore parameter uncertainty associated with predictions around this point is smaller than at the boundaries of the parameter space. As tolerance intervals incorporate parameter uncertainty, this means that the interval is smaller in the center. This is true for the factors in X as well as the mean prediction $\hat{y}$, as illustrated in Figure \ref{fig:ti-shape}. Here, observations in $X$ as well as a vector of coefficients $\beta$ were randomly generated to create the response $y$ with some added noise. After fitting an OLS model to the data, tolerance intervals were calculated for the predicted values $\hat{y}$ and their widths plotted on the y-axis. Note that this kind of curvature can be observed independently of parameter ranges or effect sizes.

While the shape can be affected by strong quadratic and interaction effects, the relationship between $\hat{y}$ and the tolerance interval width can be roughly approximated as a quadratic polynomial – a fact that is utilized in the proposed approximation. In a preprocessing step performed before DS optimization, the original CQA model is used to predict means and TI ranges for a design matrix $X$ specifically designed to detect quadratic trends in data – a central composite design (CCD) \cite{mongomery2017doe}. The design is composed of $2^p+2p+5$ rows in $X$, giving it reasonable scalability over the number of factors typically used in biopharmaceutical models. The data is then used to regress the predicted means $\hat{y}$ onto the TI widths assuming a second order polynomial:

\begin{equation}
\hat{ti}(\hat{y})=\hat{\beta}_0+\hat{\beta}_1\hat{y}+\hat{\beta}_2\hat{y}^2
\end{equation}

This regression model of TI widths is subsequently used in the optimization process. Furthermore, lower and upper boundaries are defined around the predicted mean:

\begin{align}
\hat{ti}_l(\hat{y})&=\hat{y}-\hat{ti}(\hat{y}) \\
\hat{ti}_u(\hat{y})&=\hat{y}+\hat{ti}(\hat{y})
\end{align}

\begin{figure}[H]
    \begin{subfigure}[t]{0.50\textwidth}
        \centering
        \includegraphics[width=\textwidth]{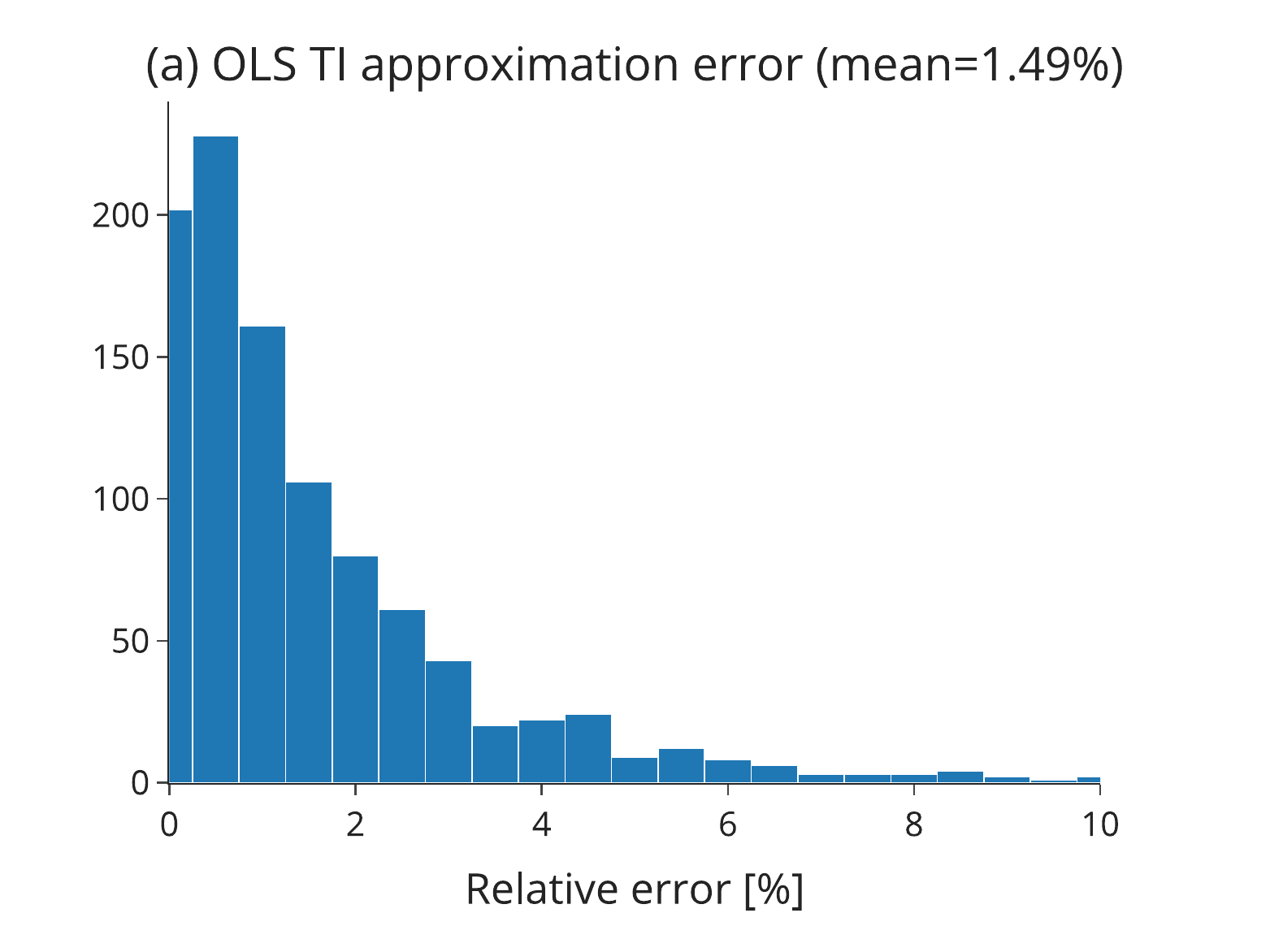}
    \end{subfigure}
    \hfill
    \begin{subfigure}[t]{0.50\textwidth}
        \centering
        \includegraphics[width=\textwidth]{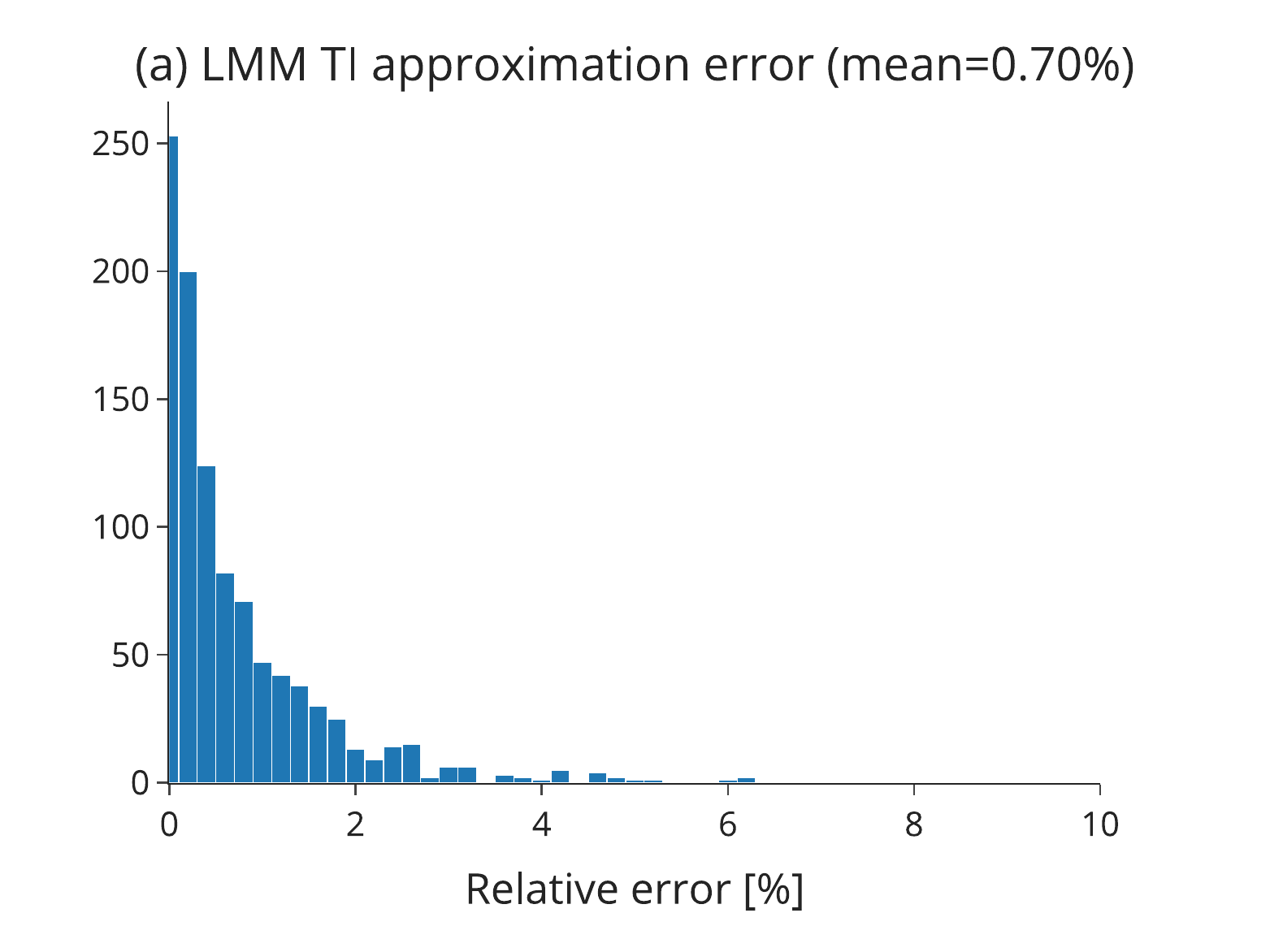}
    \end{subfigure}

\caption{Relative error of the TI approximation for (a) OLS and (B) LMM, calculated as $err=(TI_{approx}-TI_{true})/(y_{max}-y_{min})$.}
\label{fig:ti-approx-accuracy}
\end{figure}

To evaluate the accuracy of the approximation, it was compared to the actual tolerance interval over a range of randomly generated values for $p$, $\sigma^2$, $\alpha$ and $\psi$ as well as different interaction and quadratic effects. A measure of relative error was calculated by dividing the difference between real and approximated TI boundaries by the range of the model response. For OLS models, the mean error over 1000 iterations was 1.49\%, normalized over the model response range. To highlight that the approximation method is largely model agnostic, the simulation was repeated with linear mixed models (LMM) and data containing a single random effect. In each iteration the random effect variance and BLUP values were varied while the data was evenly split into four random blocks. For calculating tolerance intervals, we used the method proposed by Franzq et al. \cite{francq2019confidence}. Here, the relative approximation error is even smaller at 0.70\%, as the range of the model response used for normalization is inflated due to additional random effect variance. Figure \ref{fig:ti-approx-accuracy} contains histograms of the error distribution for both the OLS and the LMM simulation.

While this level of accuracy cannot be expected for more complex model types, we believe that it yields a reasonable estimate of the TI to guide the optimization process. Furthermore, inaccuracies due to the approximation are mitigated by the second pass of the optimizer that uses the actual TI calculation method, as described in the following sections.

\subsection{Optimization Problem and Constraints}
We define the search for a rectangular design space with the maximum volume as a continuous optimization problem with an objective function that represents the hyperrectangles volume:

\begin{equation} \label{eq:objective}
\begin{aligned}
\textrm{maximize} \quad & \prod_{i=1}^p{(x_{p+i}-x_i)w_i}, \quad & x \in \mathbb{R}^{2p}
\\
\textrm{subject to} \quad & c_1(x), ..., c_m(x) \geq 0
\end{aligned}
\end{equation}

Where $x$ is the vector of parameter values that is varied in the optimization process, comprised of the lower parameter boundaries in the first p elements and the upper boundaries in the next p elements, with p being the number of factors in the model. One can see that equation \ref{eq:objective} simply maximizes the volume of the hypercube spanned by lower and upper parameter ranges. The vector w contains predefined weights per parameter which can be used to favour one parameter range over. By default, this a vector of ones. The effects of weighting are illustrated in the example in section 4.1.

To meet all requirements for a valid design space, a total of $m=2^{p+1}+5p+2$ inequality constraints is defined. First, the parameter space to be searched is constrained by ensuring that a parameter’s lower range boundary is smaller than its upper boundary.

\begin{align}
c_i(x) &= x_{p+i}-x_i, & i=1,2,...,p
\end{align}

Oftentimes the design space is required to contain each parameter’s setpoint, so that $x_i \leq s_i \leq x_{p+1}$. This is expressed as the inequality constraint:

\begin{align}
c_i(x) &= s_i-x_i, & \quad i=1,2,...,p \\
c_i(x) &= x_i-s_i, & \quad i=p+1,...,2p
\end{align}

Similarly, optimization of the parameter space should only be performed within the screening range boundaries $b_l$ and $b_u$. As our main optimization algorithm does not support natural boundaries, this is implemented as inequality constraints:

\begin{align}
c_i(x) &= x_i-b_{l,i}, & i=1,2,...,p \\
c_i(x) &= b_{u,i}-x_i, & i=p+1,...,2p
\end{align}

The remaining constraints address the evaluation of the TI. To that end, the approximation generated in the pre-optimization step is used to calculate boundaries around CQA predictions that capture model uncertainty. These boundaries are then compared against the lower and upper acceptance limits $a_l$ and $a_u$ in each of the $2^p$ corner points of the hyperrectangle:

\begin{align}
c_i(x) &=\hat{ti}_l(x\hat{\beta})-a_l, & i=1,2,...,2^p \label{eq:ti-corner-l} \\ 
c_i(x) &= a_u-\hat{ti}_u(x\hat{\beta}), & i=1,2,...,2^p \label{eq:ti-corner-u}
\end{align}

Here, $x$ denotes a point that is taken from all corners of the current DS candidate and $x\hat{\beta}$ yields the predicted mean that is passed to the TI approximation.

Evaluating the TI at corners alone does not guarantee a valid design space, as curvature in the response surface might lead to parameter ranges between corner points that exceed acceptance limits. To resolve this problem, a final inequality constraint uses a nested optimization step to find minima and maxima of the TI boundaries inside the hypercube:

\begin{align}
c_i(x) &= \min_{x \in DS} \hat{ti}_l(x\hat{\beta})-a_l \\
c_i(x) &= a_u - \max_{x \in DS} \hat{ti}_u(x\hat{\beta})
\end{align}

While these last constraints might seem to make equations \ref{eq:ti-corner-l} and \ref{eq:ti-corner-u} redundant, having both TI checks in place can improve convergence of the optimizer in certain scenarios. Furthermore, one or the other can be deactivated in practice, depending on the type of optimization problem.

\subsection{Categorical Factors}
Many different approaches to optimization problems involving categorical factors can be found literature, e.g., mixed-integer programming \cite{munoz2020global}, branch-and-bound tree searching \cite{vanaret2021global}, genetic algorithms \cite{badran2008integrating}, Bayesian methods \cite{saves2022bayesian}, etc.

Here we describe how to include categorical variables as continuous parameters of an optimization problem by exploiting how they are represented in a regression model. Specifically, we assume that such factors enter the model as sum-coded columns of the design matrix $X$ \cite{ucla2021contrast}. A categorical factor with $k$ levels is encoded into $k-1$ columns and their values indicate that the level corresponding to the observation is active (“1”), inactive (“0”) or that the last level, which is not encoded as a separate column, should be applied (“-1”). An example for this schema is given in table \ref{tab:cat-1} and \ref{tab:cat-2}.

\begin{table}[H]
\begin{minipage}{.5\linewidth}
    \caption{Original, categorical column} \label{tab:cat-1}
    \begin{tabular}{p{0.9\textwidth}}
    \toprule
    Category \\
    \midrule
    A \\
    A \\
    B \\
    B \\
    C \\
    C \\
    \bottomrule
    \end{tabular}
\end{minipage}
\begin{minipage}{.5\linewidth}
    \caption{Sum-coded column} \label{tab:cat-2}
    \begin{tabular}{p{0.45\textwidth}p{0.45\textwidth}}
    \toprule
    Category | A & Category | B \\
    \midrule
    1 & 0 \\
    1 & 0 \\
    0 & 1 \\
    0 & 1 \\
    -1 & -1 \\
    -1 & -1 \\
    \bottomrule
    \end{tabular}
\end{minipage} 
\end{table}

\begin{figure}[H]
\centering
\includegraphics[width=0.9\textwidth]{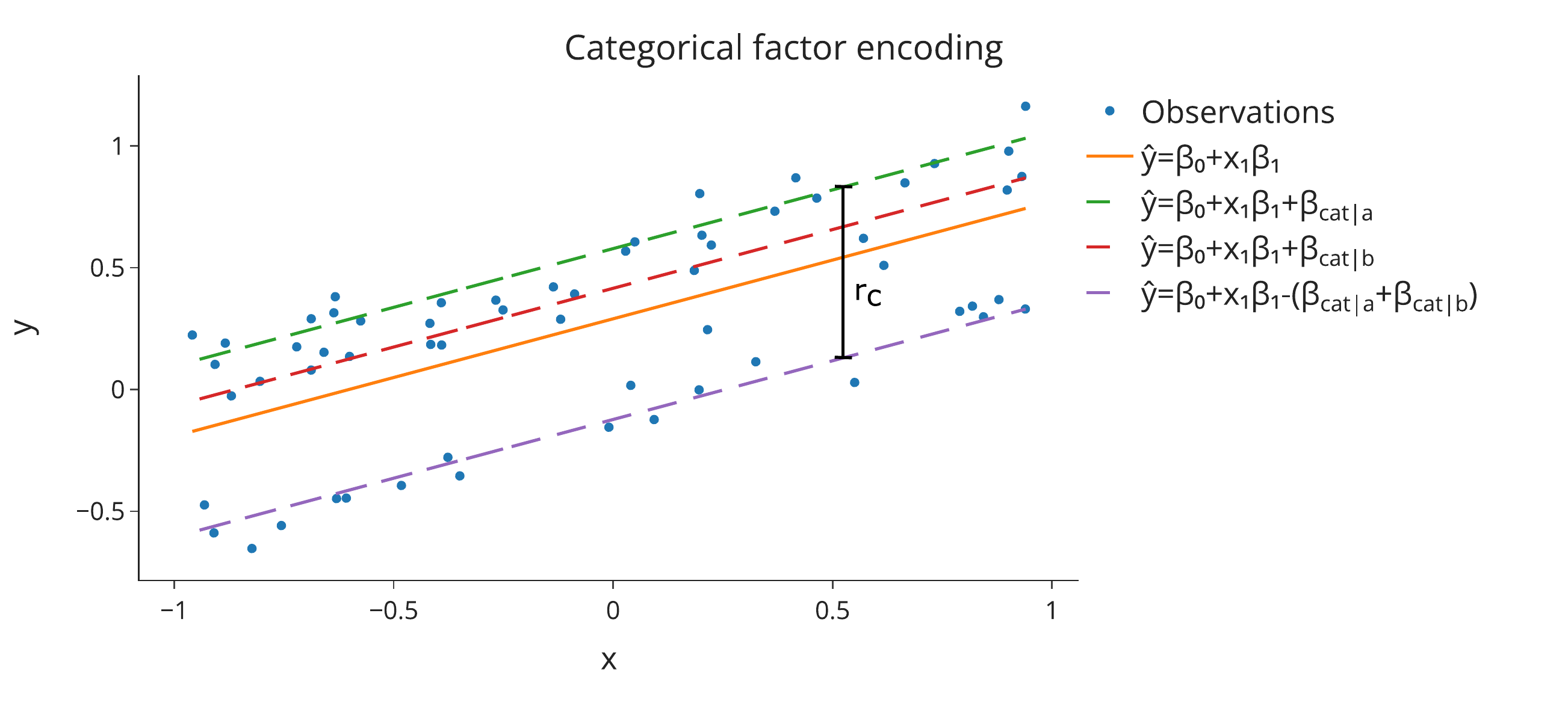}
\caption{The effect of different levels of a sum-coded, categorical factor.}
\label{fig:encoding}
\end{figure}

When representing the training data for the model in this way, the least-squares solution yields coefficients for the first $k-1$ levels of the categorical factor, hereafter denoted as the vector $\hat{\beta}_c \in \mathbb{R}^{k-1}$. The coefficient for the last level, not represented as a column in the training data, can be calculated by $-\sum_{j=1}^{k-1}\hat{\beta}_{c,j}$. Consequently, the coefficients sum to zero and their values correspond to offsets from the mean of level means, or the intercept in the case of OLS models. In other words, the effect of a particular level in a categorical factor is simply added to the model response $\hat{y}=X\hat{\beta}$. Figure \ref{fig:encoding} shows the vertical shift of the regression line caused by three categorical levels “a”, “b” and “c”, whereas “c” is encoded as the negative sum of all other level effects.

This description of a categorical factor allows us to incorporate it into the optimization problem as a continuous factor bounded by $[min(\hat{\beta}_c),max(\hat{\beta}_c)]$. The range that should be considered by the optimizer is denoted as $r_c=max(\hat{\beta}_c)-min(\hat{\beta}_c)$, shown in Figure \ref{fig:encoding}. As the optimizer operates in a normalized space, i.e., parameters are deviated within $[-1,1]$, a mean-normalization rescaling scheme is incorporated into the evaluation of $\hat{y}$ within the optimization procedure. To that end the mean prediction is expanded to $\hat{y}=X\hat{\beta}+x_c r_c + \bar{r_c}$ or, equivalently, by simply adding the terms model data and parameters $\hat{y}=[x_1,...,x_p,x_c,\textbf{1}_n]^T[\hat{\beta}_1,...,\hat{\beta}_p,\bar{r_c}]$. After optimization, a lower and upper range is returned that encloses all valid levels. These levels are added to the design space as possible values for the categorical factor.

While this approach is specific to polynomial regression models and sum-coded categorical factors, it poses fewer restrictions on the choice of optimization algorithm when compared to other methods.

%% file: sections/04_sim_results.tex
\section{Simulation Study Results}
\subsection{Accuracy}
\begin{figure}[H]
    \begin{subfigure}[t]{0.50\textwidth}
        \centering
        \includegraphics[width=\textwidth]{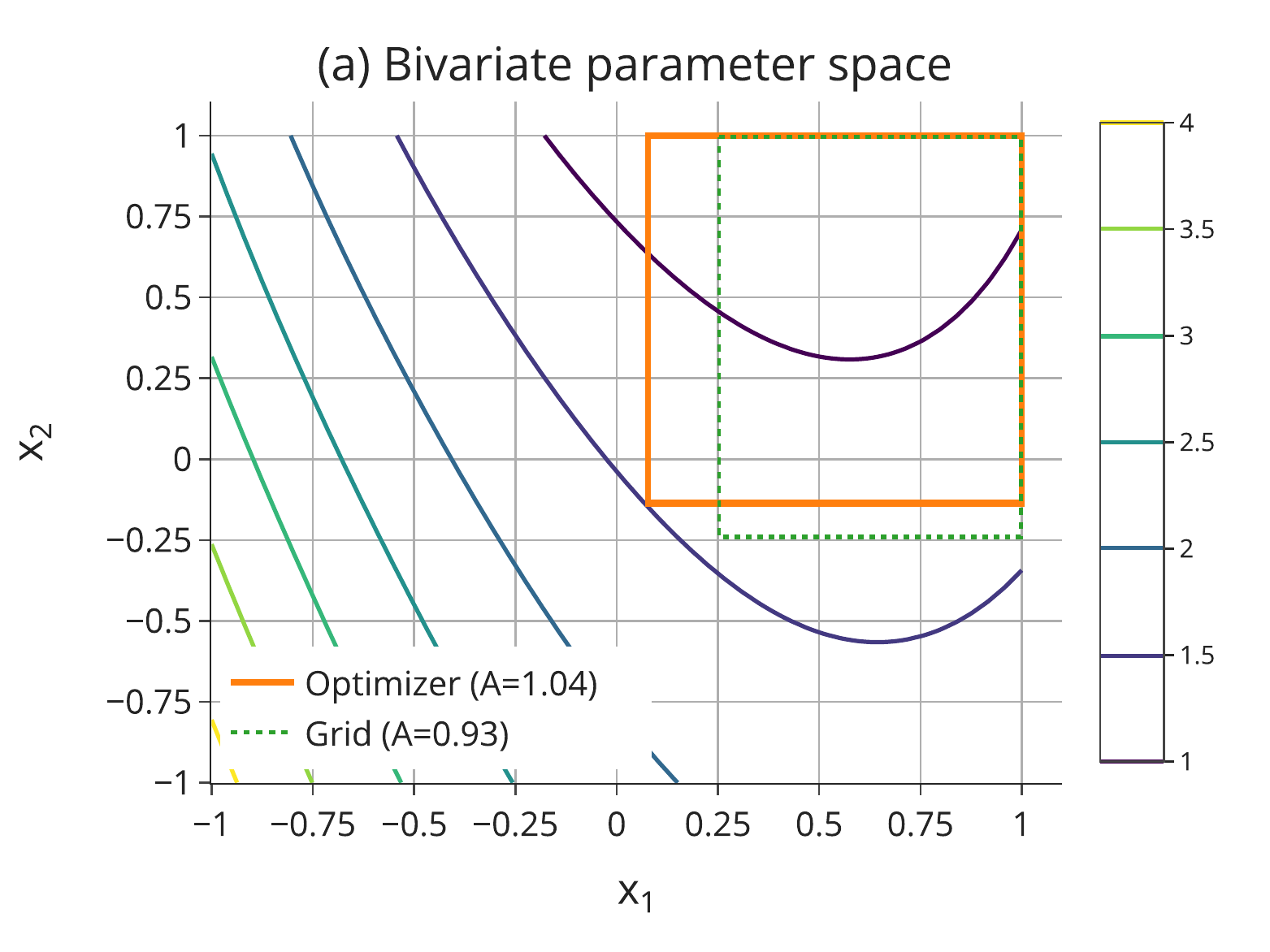}
    \end{subfigure}
    \hfill
    \begin{subfigure}[t]{0.50\textwidth}
        \centering
        \includegraphics[width=\textwidth]{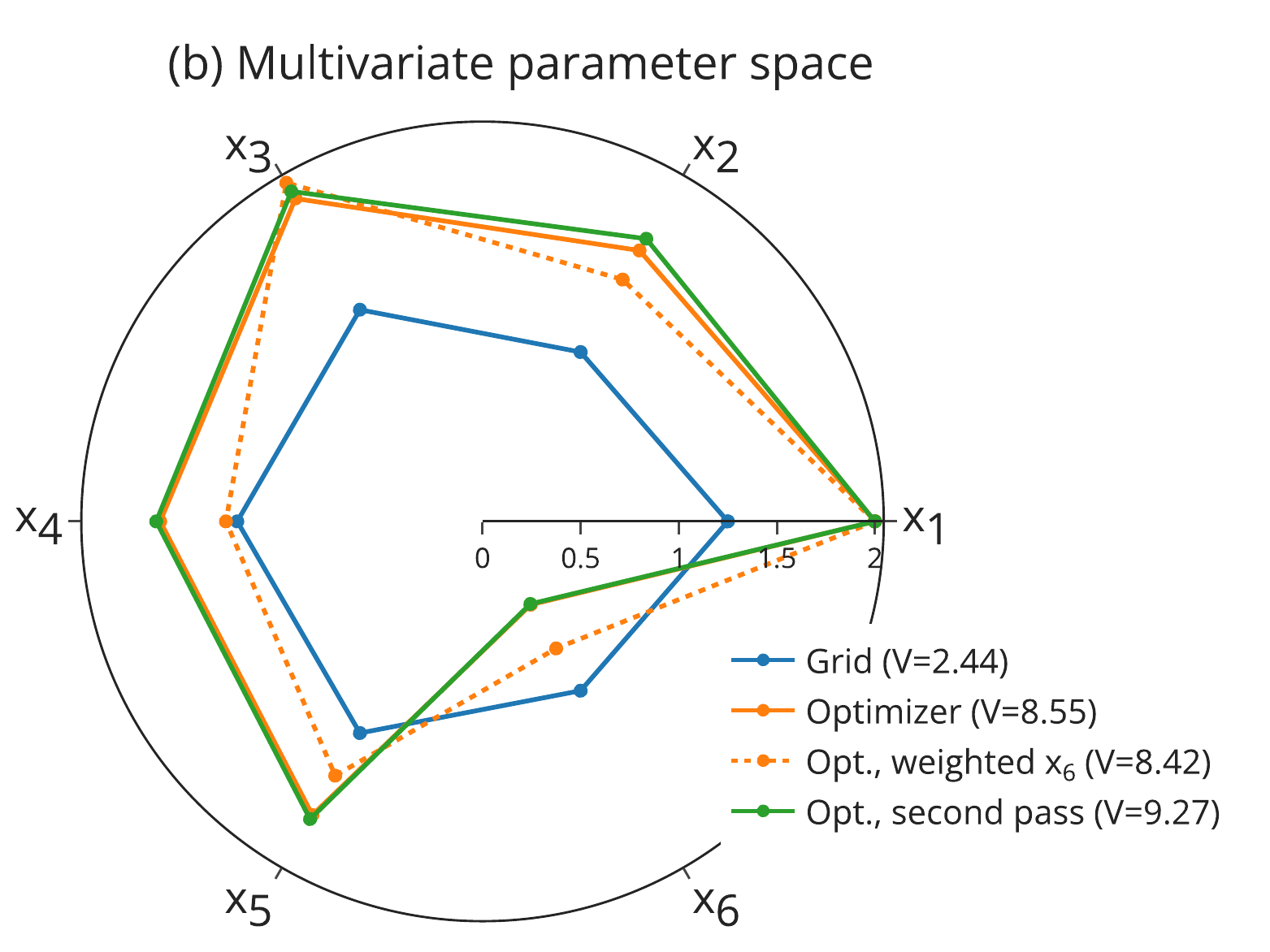}
    \end{subfigure}

\caption{(a) Bivariate parameter space with an upper tolerance interval boundary shown as the contour. The design spaces found by the optimizer and the grid method are marked as rectangles. (b) A multivariate parameter space illustrated as a spider plot with ranges per parameter shown in the axes and DS plotted as polygons.}
\label{fig:res-accuracy}
\end{figure}

The accuracy of the proposed method was evaluated in simulations using different CQA models with varying number of effects and randomized effect strengths, as well as different interactions and quadratic effects. We pick two representative examples from the results and compare them to a grid-based approach commonly found state-of-the-art software. Such DS calculation schemes discretize the parameter space by dividing it into a grid of equidistant points and evaluating all possible point combinations to find the largest hyperrectangle that contains only points satisfying the acceptance criteria. The downsides of this approach are twofold: one is a lack of accuracy as the DS is restricted to points on the grid and another concerns the computation time of such brute-force methods, explained in the next section. The accuracy problem is illustrated in Figure \ref{fig:res-accuracy}a. As the grid-based method can only evaluate points on the 8x8 grid, the resulting DS volume is smaller than that found by the optimizer, which is not subject to that limitation. Of course, the accuracy of the grid method can be improved by increasing the resolution. This, however, can drastically increase the computation time and it is not clear which grid resolution yields a DS with reasonable accuracy.

Going beyond the trivial example of two parameters, Figure 6b shows the results for a more complex model comprised of six main factors and several interaction and quadratic effects:

\begin{equation}
\begin{split}
y = -0.45+2.1x_1-0.93x_2+0.63x_3+0.42x_4-0.32x_5+0.28x_6\\
+0.76x_2^2-0.21x_1 x_2-0.34x_1 x_5+0.22x_1 x_6+0.27x_2 x_3\\
+0.24x_2 x_4-0.25x_2 x_5-0.19x_3 x_6-0.32x_4 x_6
\end{split}
\end{equation}

In this representation of the DS, valid parameter ranges are indicated by the axes of the spider plot. As in the previous example, the volume of the DS based on the $9^6$ grid is much smaller than that of the optimizer (orange solid line in Figure \ref{fig:res-accuracy}b). It is, however, more balanced, giving each parameter a similar range. This is an effect that might be observed in optimizer-based results for higher-dimensional parameter spaces, as the optimization objective only considers DS volume. Here, the range of $x_6$ is quite small compared to the grid method. In case such a result is not desired, the weight for $x_6$ in the objective equation \ref{eq:objective} can be increased to yield a more balanced parameter range while sacrificing some volume, indicated by the dashed orange line in Figure \ref{fig:res-accuracy}b. Both DS plotted in orange were calculated using the first phase of the proposed scheme only, i.e., with an approximated TI and $\rho_{start}=0.01$. The unweighted DS was then passed to the second phase for further refinement using the true TIs and $\rho_{start}=0.001$, which results in the largest volume. The parameter ranges of the DS are shown in green and numeric results are summarized in table \ref{tab:res-accuracy}.

\begin{table}[H]
\caption{Accuracy evaluation results} \label{tab:res-accuracy}
\small
{\renewcommand{\arraystretch}{1.3}
\begin{tabular}{llllllllllllll}
\hline
                           & \multicolumn{2}{c}{\textbf{x1}} & \multicolumn{2}{c}{\textbf{x2}} & \multicolumn{2}{c}{\textbf{x3}} & \multicolumn{2}{c}{\textbf{x4}} & \multicolumn{2}{c}{\textbf{x5}} & \multicolumn{2}{c}{\textbf{x6}} & \textbf{volume} \\
                           & \textbf{from}   & \textbf{to}   & \textbf{from}   & \textbf{to}   & \textbf{from}   & \textbf{to}   & \textbf{from}   & \textbf{to}   & \textbf{from}   & \textbf{to}   & \textbf{from}   & \textbf{to}   &                 \\ \hline
\textbf{grid}              & -1.00           & 0.00          & 0.00            & 1.00          & -1.00           & 0.25          & -0.75           & 0.50          & -0.25           & 1.00          & -0.25           & 1.00          & 2.44            \\
\textbf{optimizer}         & -1.00           & -0.51         & -0.63           & 0.97          & -1.00           & 0.64          & -1.00           & 0.73          & -0.91           & 1.00          & -1.00           & 1.00          & 8.55            \\
\textbf{opt., weighted}    & -1.00           & -0.51         & -0.69           & 0.98          & -1.00           & 0.66          & -1.00           & 0.76          & -0.95           & 1.00          & -1.00           & 1.00          & 8.42            \\
\textbf{opt., second pass} & -1.00           & -0.25         & -0.43           & 1.00          & -1.00           & 0.31          & -1.00           & 0.50          & -1.00           & 1.00          & -1.00           & 1.00          & 9.27            \\ \hline
\end{tabular}}
\end{table}

\subsection{Computing Time}
\begin{figure}[H]
\centering
\includegraphics[width=0.9\textwidth]{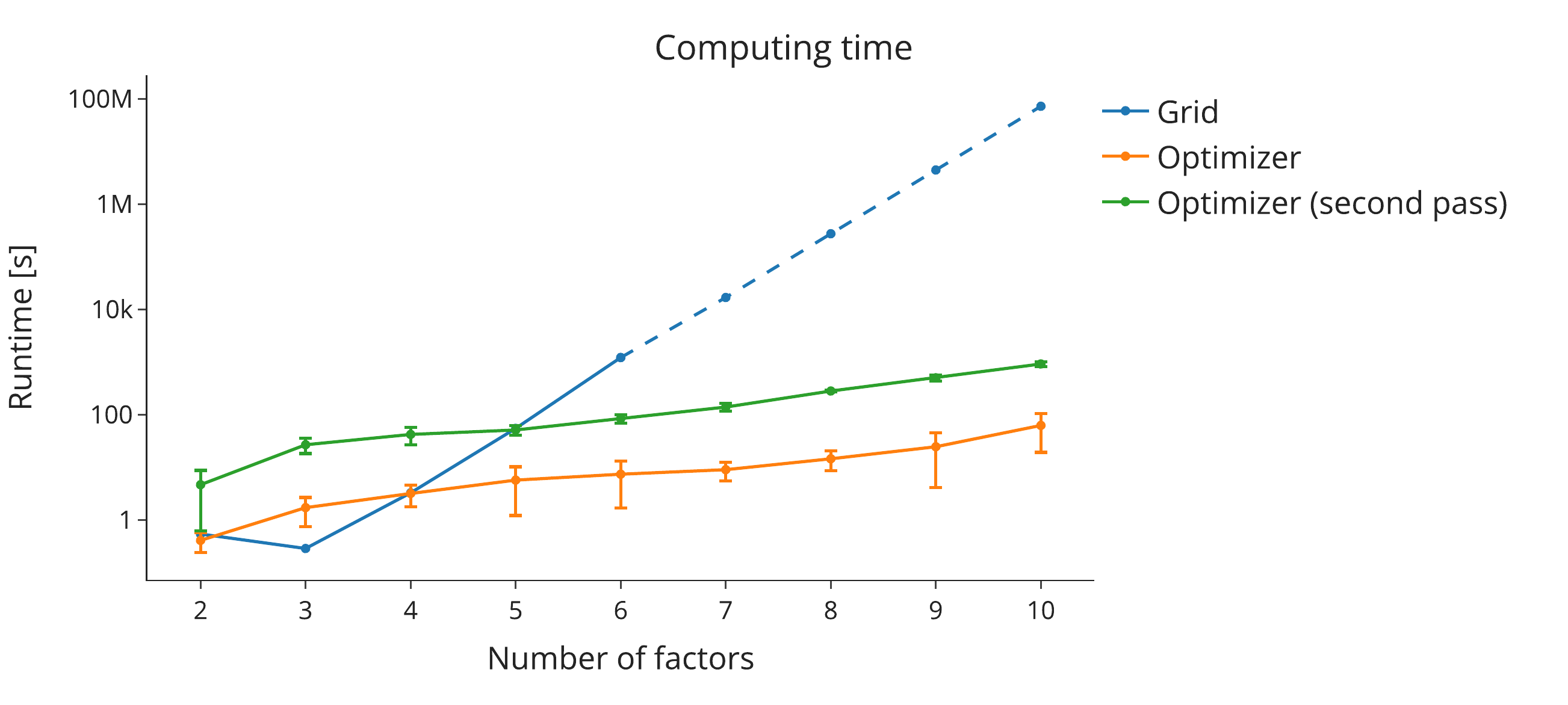}
\caption{Performance evaluation of the proposed method using both the TI approximation and the true TI, as well as the reference method.}
\label{fig:performance}
\end{figure}

A major problem of grid-based approaches is poor scaling to higher dimensions of the parameter space. As the space to be searched grows exponentially with the number of parameters, so does computing time required to find the DS. This renders the method infeasible for many use cases in the biopharmaceutical domain where ten main factors or more are not uncommon for unit operation models. For example, the screening algorithm would need to evaluate $8^10$ points for parameter ranges segmented into eight parts and ten parameters while also computing all possible hyperrectangles within this space. Although optimizer-based methods are not guaranteed to find the global maximum, the parameter space is screened in a more systematic way, ideally eliminating large parts of it in a few iterations. This is illustrated in Figure \ref{fig:performance}, where the time for computing a DS is plotted against the number of model parameters of randomly generated models. For the optimizer-based method, 100 OLS models were generated per model size with different effect sizes and interactions. The plot shows the mean time over those 100 iterations as well as the standard deviation as an error bar. Due to the performance limitations described above, this was not possible for the grid-based method and only the results from one exemplary model sequence are plotted. For the same reason, those results were computed for up to six parameters and an optimistic projection is used for parameter numbers seven to ten, indicated by the dashed line.

The number of iterations for the second optimization phase was capped at 100 and the computing time added to the calculation scales linearly with the number of parameters, shown as the green line in the plot. Overall, the optimizer-based computation for 10 factors is six to seven orders of magnitude faster than the grid method, effectively enabling the exact computation of a DS in higher dimensions.

%% file: sections/05_discussion.tex
\section{Discussion}
We believe that the optimizer-based approach to calculating a DS circumvents the problems associated with grid-based methods, i.e., poor scaling over the number of parameters and a lack of accuracy due to the limitation to grid points in the solution space. In contrast to existing Bayesian methods  \cite{kusumo2019bayesian} the resulting DS is exact, confidence-free and consists of a linear combination of parameter ranges, facilitating operational simplicity for biopharmaceutical control strategies. However, as it is the case for many optimization problems - especially in high-dimensional parameter spaces – there is a possibility of converging to local minima and not finding the largest possible DS. This problem can be mitigated by using different random seeds, different starting points for the algorithm or different parameter weights. Results from different optimizer starting conditions could be compared to see if they agree and thus gain confidence in the result or only the maximum DS from the set of candidates could be returned.

As illustrated in figure \ref{fig:res-accuracy}, the weighting scheme in combination with fast computing times also facilitates the iterative exploration of the DS in higher dimensions. One can investigate the effect of weighing a parameter on the possible ranges of other parameters and thereby gain understanding about DS-related dependencies that cannot be directly derived from the model.

%% file: sections/06_conclusion.tex
\section{Conclusion}
In this contribution, we outlined a novel method for finding a ICH Q8 compliant design space comprised of linear combinations of PP ranges. The relationship between CQA and PPs is represented as a polynomial regression model and its prediction is used to evaluate whether CQAs meet acceptance criteria.

Conservative estimation methods are vital in the biopharmaceutical domain, which is why the boundaries of tolerance intervals are used for evaluation of ACs instead of the predicted mean CQA. As the TI calculation can be complex, an approximation is generated in a pre-optimization step that can be used to calculate the TI by performing a single matrix multiplication. COBYLA is used for the minimization of the main objective, and we suggest refining results with SLSQP or COBYLA with smaller $\rho_{start}$ and using the true TI instead of an approximation. 

Performance evaluations show that the proposed method results in design spaces with a larger volume when compared to PP space discretization methods, and that they can be calculated in a fraction of the time. This, for the first time, enables the calculation of design spaces for more than 10 process parameters of one model. We believe that this approach will facilitate a robust definition of the design space for biopharmaceutical development that reduces patient risk by employing conservative estimators while allowing manufacturers to maximize control ranges and operational flexibility. Furthermore, this increased flexibility, gained solely by an improved evaluation of existing data and models is a change for manufacturers to optimize a process within the DS.

%% file: sections/07_acknowledgements.tex
\section{Acknowledgements}
This work was conducted within the COMET Centre CHASE, funded within the COMET - Competence Centers for Excellent Technologies program by the BMK, the BMDW and the Federal Provinces of Upper Austria and Vienna. The COMET program is managed by the Austrian Research Promotion Agency (FFG). 

The authors acknowledge TU Wien Bibliothek for financial support through its Open Access Funding Program.